\documentclass[11pt]{article}
\usepackage[a4paper]{anysize}\marginsize{3.5cm}{3.5cm}{1.3cm}{2cm}
\pdfpagewidth=\paperwidth \pdfpageheight=\paperheight
\usepackage{amsfonts,amssymb,amsthm,amsmath,eucal}
\usepackage{pgf}
\usepackage{bbm}
\theoremstyle{plain}
\newtheorem{thm}{Theorem}[section]
\newtheorem{theorem}[thm]{Theorem}

\newtheorem{lemma}[thm]{Lemma}

\newtheorem{corollary}[thm]{Corollary}
\newtheorem{conjecture}[thm]{Conjecture}

\theoremstyle{definition}

\newtheorem{remark}[thm]{Remark}
\newtheorem{example}[thm]{Example}

\newtheorem{thevarthm}[thm]{\varthmname}

\newenvironment{varthm*}[1]{\trivlist\item[]{\bf #1.}\it}{\endtrivlist}

\renewcommand\geq{\geqslant}

\renewcommand\leq{\leqslant}

\newcommand\be{\begin{eqnarray*}}
\newcommand\ee{\end{eqnarray*}}

\newcommand\Q{\mathbb Q}

\newcommand\Z{\mathbb Z}

\newcommand\K{\mathbb K}
\renewcommand\P{\mathbb P}

\newcommand\calo{{\mathcal O}}
\newcommand\cali{{\mathcal I}}

\newcommand\newop[2]{\def#1{\mathop{\rm #2}\nolimits}}
\newop\log{log}
\newop\ord{ord}
\newop\Gal{Gal}
\newop\SL{SL}
\newop\Bl{Bl}
\newop\mult{mult}
\newop\mass{mass}
\newop\div{div}
\newop\codim{codim}
\newop\sing{sing}
\newcommand\eqnref[1]{(\ref{#1})}

\newcommand\subm{\underline{m}}
\newcommand\subn{\underline{n}}
\newcommand\one{\mathbbm{1}}
\newcommand\ldk{(\P^2,\calo_{\P^2}(d)\otimes\cali(kZ))}
\newcommand\ldmi{(\P^2,\calo_{\P^2}(d)\otimes\cali(\subm Z))}
\newcommand\ldthreem{(\P^2,\calo_{\P^2}(\widetilde{d}-3)\otimes\cali((\widetilde{\subm}-\one)Z)}
\newcommand\ldmij{(\P^2,\calo_{\P^2}(d-1)\otimes\cali((\subm-\one)Z)}
\newcommand\alphasubmn{\alpha_{\subm,\subn}}
\newcommand\mj{m^{(j)}}
\newcommand\nj{n^{(j)}}
\newcommand\mji{\mj_i}
\newcommand\submj{\subm^{(j)}}
\newcommand\subnj{\subn^{(j)}}

\begin{document}

\author{M.~Dumnicki, T.~Szemberg\footnote{The second named author was partially supported
by NCN grant UMO-2011/01/B/ST1/04875}, H.~Tutaj-Gasi\'nska}
\title{Symbolic powers of planar point configurations}
\date{\today}
\maketitle
\thispagestyle{empty}

\begin{abstract}
   Let $Z$ be a finite set of points in the projective plane and let $I=I(Z)$ be its
   homogeneous ideal. In this note we study the sequence $\alpha(I^{(m)})$, $m=1,2,3\dots,$
   of initial degrees of symbolic powers of $I$. We show how bounds on the
   growing order of elements in this sequence determine the geometry of $Z$.
\end{abstract}


\section{Introduction}
   There has been considerable interest in symbolic powers of ideals
   motivated by various problems in Algebraic Geometry, Commutative Algebra
   and Combinatorics, see e.g. \cite{ELS01}, \cite{HH07}, \cite{Sul08}
   and references therein. Our motivation originates in the theory
   of linear systems, more exactly in
   postulation or interpolation problems.
   The unsolved up to date conjectures of Nagata and Segre-Harbourne-Gimigliano-Hirschowitz
   are prominent examples of such problems. We refer to
   \cite{RDLS} for precise formulation, motivation and more open problems.
   These problems can be formulated in the language of symbolic powers.
   For example, the Nagata conjecture predicts that if $Z$ is a set
   of $r\geq 10$ \emph{general} points in the projective plane $\P^2$ and
   $I$ is the saturated ideal of $Z$, then the initial degree $\alpha(I^{(m)})$
   of the $m$--th symbolic power of $I$ is bounded from below by $m\sqrt{r}$.

   In the present paper we look at a fixed finite set $Z$ of \emph{arbitrary} points in
   the projective plane $\P^2$ with defining saturated ideal $I=I(Z)$ and we consider the increasing
   sequence
   $$\alpha(I)<\alpha(I^{(2)})<\alpha(I^{(3)})<\dots $$
   of initial degrees of symbolic powers of $I$.
   The differences in this sequence will be denoted by
   \begin{equation}\label{eq:def of difference}
      \alpha_{m,n}(Z):=\alpha(I^{(m)})-\alpha(I^{(n)})
   \end{equation}
   for $m>n$.
   Bocci and Chiantini initiated in \cite{BocCha11} the study of the relationship
   between
   the numbers $\alpha_{m,n}(Z)$ and the geometry of the set $Z$ by studying the difference $\alpha_{2,1}(Z)$.
   In particular,
   they proved that if this number is the minimal possible i.e.
   $\alpha_{2,1}(Z)=1$, then the points in $Z$ are either collinear
   or they are all the intersection points of $\alpha(2Z)$ general lines.
   It is natural to look at this problem more generally and consider also
   other numbers $\alpha_{m,n}(Z)$. Our main result is the following characterization
\begin{varthm*}{Theorem}
   If
   $$\alpha_{2,1}(Z)=\dots=\alpha_{k+1,k}(Z)=d,$$
   then
   \begin{itemize}
   \item[i)] for $d=1$ and $k\geq 2$ the set $Z$ is contained in a line, i.e. $\alpha(Z)=1$;
   \item[ii)] for $d=2$ and $k\geq 4$ the set $Z$ is contained in a conic, i.e. $\alpha(Z)=2$.
   \end{itemize}
   Moreover, both results are sharp, i.e. there are examples showing that
   one cannot relax the assumptions on $k$.
\end{varthm*}
   This result is proved in the text in  Theorems \ref{thm:first thm} and \ref{thm:last thm}.
\begin{varthm*}{Remark}
   \rm
   Naively one could expect that assuming sufficiently many equalities
   \begin{equation}\label{eq:equal alphas}
      \alpha_{2,1}(Z)=\dots=\alpha_{k+1,k}(Z)=d,
   \end{equation}
   one could always conclude $\alpha(Z)=d$. Example \ref{ex:Nagata}
   shows that this is false for all $d\geq 4$ even if the equalities
   in \eqnref{eq:equal alphas} hold for \emph{all} $k\geq 1$.

   We expect that the Theorem holds also for $d=3$ but it seems that a proof
   requires some new ideas. We hope to come back to this question
   in the future.
\end{varthm*}
   We work over an algebraically closed field of characteristic zero.

\section{Symbolic powers and initial degrees}
   Let $Z=\left\{P_1,\dots,P_r\right\}$
   be a fixed finite set of points in the projective plane. Let $I(P_i)$ be the maximal ideal
   of the point $P_i$ in the graded ring $R=\K[x_0,x_1,x_2]$, where $\K$ is an algebraically
   closed field of characteristic zero. Then the ideal of $Z$ is
   $$I=I(Z)=I(P_1)\cap\dots\cap I(P_r).$$
   For a positive integer $k$, the $k$--th symbolic power of $I$ is defined by
   \begin{equation}\label{eq:symb power}
      I^{(k)}:=I(P_1)^k\cap\dots\cap I(P_r)^k.
   \end{equation}
   The definition of the symbolic power of an arbitrary ideal is more
   complicated and involves associated primes, see \cite[Chapter IV.12, Definition]{ZarSam75},
   but in our context \eqnref{eq:symb power} suffices.
   The subscheme of $\P^2$ defined by the ideal $I^{(k)}$ will be denoted by $kZ$.

   It is convenient to study slightly more general subschemes, the so called fat points.
   To this end let $\subm=(m_1,\dots,m_r)\in\Z^r$ be a vector of non-negative
   integers. We write
   $$I(\subm Z)=I(m_1P_1+\dots+m_rP_r)=I(P_1)^{m_1}\cap\dots\cap I(P_r)^{m_r}.$$

   For an arbitrary homogeneous ideal
   $$I=\bigoplus\limits_{n\geq 0} I_n,$$
   we define the \emph{initial degree} of $I$ to be the number
   $$\alpha(I)=\min\left\{n\in \Z:\; I_n\neq 0\right\}.$$
   If $C$ is a divisor defined by some polynomial $f\in I_n$, then we
   abuse a little bit the notation and write $C\in I_n$. As the explicit
   equations defining a divisor will never play a role in this paper, this should cause no confusion.

   Finally, we write $\one$ for the vector $(1,\dots,1)\in\Z^r$ and we extend
   the definition in \eqnref{eq:def of difference} to the case of inhomogeneous
   multiplicities by putting
   $$\alphasubmn(Z):=\alpha(I(\subm Z))-\alpha(I(\subn Z)).$$
   Thus $\alpha_{m,n}(Z)=\alphasubmn(Z)$ for $\subm=(m,\dots,m)$ and $\subn=(n,\dots,n)\in\Z^r$.

   We begin by the following observation which is crucial in the sequel. Roughly speaking,
   it amounts to saying that if a divisor $C$ is of minimal degree and some difference
   $\alpha$ is given, then every effective subdivisor of $C$, in particular every irreducible component of $C$,
   is again of minimal degree
   for an appropriate subscheme and the difference at this subscheme is at most that at $Z$.
\begin{lemma}\label{lem:alphas on subdivisors}
   Let $Z$ be a fixed set of points $P_1,\dots,P_r\in\P^2$ and
   let $m>n$ be positive integers. Let
   $\beta=\alpha(I(mZ))$, $\gamma=\alpha(I(nZ))$ and
   $\alpha=\alpha_{m,n}(Z)=\beta-\gamma$. Let $C\in I(mZ)_{\beta}$
   be an effective divisor. Furthermore let
   $$C=C_1+C_2$$
   be a sum of two integral non-zero divisors.
   Let $\beta_j=\deg(C_j)$, $m_i^{(j)}=\ord_{P_i}C_j$,
   $\submj=(\mj_1,\dots,\mj_r)$,
   $$n_i^{(j)}=\max\left\{\mj_i-(m-n),\, 0\right\}$$
   and $\subnj=(\nj_1,\dots,\nj_r)$.
   Then
   \begin{itemize}
   \item[i)] $\beta_j=\alpha(I(\submj Z))$ and
   \item[ii)] $\alpha_{\submj,\subnj}\leq\alpha$.
   \end{itemize}
   for $j=1,2$.
\end{lemma}
\proof
   We prove the statement for $C_1$, the claim for $C_2$ follows by symmetry.

   i) Assume to the contrary that there is a divisor $C_1'$
   of degree $\deg(C_1')=\beta_1'<\beta_1$ with vanishing vector $\subm^{(1)}$.
   Then $C_1'+C_2$ has vanishing vector $\subm^{(1)}+\subm^{(2)}=\subm$
   and degree $\beta_1'+\beta_2\,<\,\beta_1+\beta_2\,=\,\beta$,
   a contradiction.

   ii) If $\subn^{(1)}=0$, then let $\Gamma_1=0$, otherwise let $\Gamma_1$ be a nonzero
   divisor of least degree $\gamma_1$ in $I(\subn^{(1)}Z)$. Let $D_1=\Gamma_1+C_2$.
   Then
   $$\ord_{P_i}D_1=\ord_{P_1}\Gamma_1+\ord_{P_i}C_2=\max\left\{m^{(1)}_i-(m-n),\, 0\right\}+m^{(2)}_i\geq n.$$
   Indeed, if $m^{(1)}_i-(m-n)=\max\left\{m^{(1)}_i-(m-n),\, 0\right\}$, then
   $$m^{(1)}_i-(m-n)+m^{(2)}_i=m^{(1)}_i+m^{(2)}_i-m+n\geq n.$$
   If it is $0=\max\left\{m^{(1)}_i-(m-n),\, 0\right\}$, then
   $$0\geq m^{(1)}_i-(m-n)\geq n-m^{(2)}_i$$
   implies
   $$n\leq m^{(2)}_i=\ord_{P_i}C_2=\ord_{P_i}D_1.$$
   Hence it must be
   $$\deg(D_1)=\gamma_1+\beta_2\geq \gamma$$
   by assumption, which implies
   $$\alpha_{\subm^{(1)},\subn^{(1)}}=\beta_1-\gamma_1=\beta-(\beta_2+\gamma_1)\leq \beta-\gamma=\alpha$$
   as asserted.

   Note that it might happen that some entries in $\subnj$ are zero or
   even that both vectors $\subn^{(1)}$ and
   $\subn^{(2)}$ are zero.
\endproof
\begin{corollary}\label{cor:alphas of components}
   Keeping the notation from the above Lemma, let
   $C$ be a divisor in $I(mZ)_{\beta}$ with the decomposition
   \begin{equation*} 
      C=\sum a_jC_j
   \end{equation*}
   into distinct irreducible
   components with $\gamma_j=\deg(C_j)$ and $\mji=\ord_{P_i}C_j$.
   Let $\submj=(\mj_1,\dots,\mj_r)$. Then
   \begin{itemize}
   \item[i)] $\gamma_j=\alpha(I(\submj Z))$ and
   \item[ii)]  for $n^{(j)}_i=\max\left\{\mj_i-(m-n),\, 0\right\}$ we have $\alpha_{\submj,\subnj}(Z)\leq \alpha$.
   \end{itemize}
\end{corollary}
\proof
   Simply apply Lemma \ref{lem:alphas on subdivisors} to every component $C_j$ of $C$
   and its residual divisor $C-C_j$.
\endproof

\section{Differences equal $1$}
   We begin by comparing $\alpha(kZ)$ and $\alpha(Z)$, which is
   the first natural generalization of the situation studied
   by Bocci and Chiantini. Since a derivative of a polynomial
   having $k$--fold zeroes at $Z$ has zeroes of order at least
   $k-1$, we have
\begin{equation}\label{eq:1}
   \alpha_{k,1}(Z)\geq k-1.
\end{equation}
   Now, we describe all cases when there is the equality in \eqnref{eq:1},
   the case $k=2$ being settled by Bocci and Chiantini.
\begin{theorem}\label{thm:first thm}
   Let $Z\subset\P^2$ be a finite set of points and let $k\geq 3$
   be an integer. If
   $$\alpha_{k,1}(Z)=k-1,$$
   then the points in $Z$ are collinear.
\end{theorem}
\proof
   The idea is to reduce the statement to the classification of Bocci
   and Chiantini and to exclude the configuration of intersection points of lines.

   Let $C$ be a divisor of degree $d=\alpha(kZ)$ vanishing at $kZ$, defined by
   the polynomial $f$. Let $D$ be the divisor given by some derivative
   of $f$ of order $(k-2)$. Then $D$ is a divisor of degree $d-(k-2)$
   vanishing along $2Z$. Our assumption
   $$\alpha(Z)=\alpha(kZ)-\alpha_{k,1}(Z)=d-(k-1)$$
   implies that $\alpha_{2,1}(Z)=1$. Thus \cite[Theorem 1.1]{BocCha11} applies
   and either all points in $Z$ are collinear in which case we are done, or
   they are intersection points of $(d-(k-2))$ general lines. In the later case, let
   $L$ be one of these lines. Then $L\cdot C=d$ and
   \begin{itemize}
   \item[a)] either $L\cdot C\geq (d-(k-2)-1)k$,
   \item[b)] or $L$ is a component of $C$.
   \end{itemize}
   In case a) we see immediately that
   $$k\geq d.$$
   Since $C$ contains points of multiplicity $k$, it follows, that in fact
   $k=d$. But then $Z$ consists of a single point and we are done.

   So we may assume that the case b) holds for all lines in the arrangement.
   Subtracting all these lines from $C$, we obtain a curve $C'$ of degree
   $d'=d-(d-(k-2))=k-2$ passing through $(k-2)Z$. Hence $d'\geq k-2$ which gives again
   the multiplicity equal to the degree and we are done.
\endproof
\begin{remark}
   One could revoke a result of W\"ustholz \cite[page 76]{EV83} in order to establish
   the above theorem. We preferred however to present an elementary proof.
   Note also for the future reference that the results along the lines of
   \cite{EV83} do not lead to a proof of Theorem \ref{thm:last thm}.
\end{remark}
   Another possible approach to a generalization of the result of Bocci and Chiantini
   is to replace the number $\alpha_{2,1}$ by $\alpha_{k,k-1}$. This leads to some
   new configurations.

   We begin by the following general Lemma.
\begin{lemma}\label{lem:only line}
   Let $Z=\left\{\,P_1,\dots,P_r\,\right\}$ be a set of $r$ points in the projective plane $\P^2$
   and let $m_1,\dots,m_r$ be positive integers. Suppose that the minimal
   degree of a divisor passing through points $P_1,\dots,P_r$ with
   multiplicities at least $m_1,\dots,m_r$ is $d$, and for multiplicities
   $(m_1-1),\dots,(m_r-1)$, it is $d-1$. Moreover we assume that there
   is a reduced and irreducible curve
   $$C\in H^0\ldmi.$$
   Then $C$ must be a line (and consequently $d=m_1=\dots=m_r=1$).
\end{lemma}
\proof
   Assume to the contrary that the degree $d$ of $C$ is at least $2$.
   Let $C'$ be a divisor defined by a first order derivative
   of the equation of $C$. Then
   $$C'\in H^0\ldmij$$
   and $C$ and $C'$ have no common
   components.
   It follows that
   \begin{equation}\label{eq:7}
      d(d-1)=C\cdot C'\geq \sum\limits_{i=1}^rm_i(m_i-1).
   \end{equation}
   On the other hand, there is by assumption no curve of degree $(d-2)$
   passing through $P_1,\dots,P_r$ with multiplicities
   $(m_1-1),\dots,(m_r-1)$. A simple dimension count gives
   \begin{equation}\label{eq:8}
      d(d-1)\leq \sum\limits_{i=1}^rm_i(m_i-1).
   \end{equation}
   It follows that
   \begin{equation}\label{eq:9}
      d(d-1)=\sum\limits_{i=1}^rm_i(m_i-1).
   \end{equation}
   In particular this means that the divisors $C$ and $C'$ meet only in points $P_1,\dots,P_r$.

   Note that \eqnref{eq:9} implies that there are at least $d$ independent sections
   in $H^0\ldmij$. Indeed,
   $$h^0\ldmij\geq {{d+1}\choose{2}}-\sum\limits_{i=1}^r{m_i\choose2}=d.$$
   Hence there is a section $s''\in H^0\ldmij$ vanishing at a point $R\in C$
   different from $P_1,\dots,P_r$. By \eqnref{eq:9} the divisor $C''$ defined by $s''$
   must have a common component with $C$. As $C$ is irreducible and the degrees
   do not agree, this is clearly impossible.
\endproof
   The key idea in the proof of the next theorem is to observe that
   the above Lemma applies to reducible curves component by component.
\begin{theorem}\label{thm:only lines}
   Let $Z=\left\{P_1,\dots,P_r\right\}\subset\P^2$ be a finite set of points. If
   $$\alpha_{k,k-1}(Z)=1$$
   for some $k\geq 2$ then
   \begin{itemize}
   \item[a)] either the points in $Z$ are collinear,
   \item[b)] or
   $Z$ consists of all intersection points of some arrangement of lines
   (which might be non reduced and degenerate i.e. containing points
   through which more than two lines pass).
   \end{itemize}
\end{theorem}
\proof
   Let $C$ be a divisor of minimal degree $d$ vanishing at $kZ$
   and let
   $$C=\sum\limits_{j=1}^s a_jC_j$$
   be its decomposition with irreducible and reduced components $C_j$.
   Let $\gamma_j$ be the degree of $C_j$ and
   $m^{(j)}_i=\ord_{P_i}C_j$ its order of vanishing at the point $P_i$.

   It follows from Corollary \ref{cor:alphas of components}
   that $\gamma_j$ is the least number $\gamma$ such that there exists
   a curve $D$ of degree $\gamma$ vanishing at $P_1,\dots,P_r$ with
   multiplicities $m^{(j)}_1,\dots,m^{(j)}_r$.

   The same Corollary implies that $\gamma_j-1$ is the least number $\gamma$ such that there exists
   a curve $D$ of degree $\gamma$ vanishing at $P_1,\dots,P_r$ with
   multiplicities $(m^{(j)}_1-1),\dots,(m^{(j)}_r-1)$.

   Thus we can apply Lemma \ref{lem:only line} to every component $C_j$
   and we are done as we get case a) for $s=1$ and case b) otherwise.
\endproof
   We conclude this section with the following corollary, strengthening
   the statement of Theorem \ref{thm:first thm}.
\begin{corollary}\label{cor:collinear}
   Assume that $k\geq 3$ and $Z$ is a collection of points such that
   $$\alpha_{k,k-1}(Z)\,=\, \alpha_{k-1,k-2}(Z)=1$$
   Then $Z$ consists of collinear points.
\end{corollary}
\proof
   Let $C\in I(kZ)$ be a divisor of degree $d=\alpha(kZ)$.
   By Theorem \ref{thm:only lines} $C$ consists of $d$ lines.
   If among these lines there are $3$ (or more) in general position,
   i.e. not all passing through the same point, then we can
   subtract these $3$ lines from $C$ obtaining a curve vanishing
   along $(k-2)Z$ of degree $d-3$, which contradicts our assumptions.
   Thus all lines in $C$ must pass through the same point. Assume
   that $C_{red}$ the support of $C$ consists of $a$ lines.
   Subtracting all of them once from $C$ produces a curve of degree $d-a$
   vanishing along $(k-1)Z$. This implies $a=1$ and we are done.
\endproof
   Note that the following example shows that it is essential that
   there are two \emph{consecutive} differences of $1$ in the sequence $\alpha(kZ)$.
\begin{example}
   Let $Z=\left\{P_1,P_2,P_3\right\}$ be $3$ general points.
   Then
   $$\alpha(kZ)=\left\{\begin{array}{ccl}
      3m-1 & \mbox{ for } & k=2m-1 \mbox{ odd}\\
      3m   & \mbox{ for } & k=2m \mbox{ even}
      \end{array}\right.$$
   In particular infinitely many differences $\alpha_{k,k-1}$ are equal $1$.
\end{example}
\proof
   This is obvious for $m=1$ and we proceed by induction on $m$.
   Assuming the theorem for $2m-1$ and $2m$, we want to prove it
   for $2(m+1)-1=2m+1$ and $2(m+1)$. In the even case we apply
   Bezout theorem in order to show that a divisor vanishing
   along $Z$ to order $2(m+1)$ must be the $(m+1)$--fold
   union of the $3$ lines determined by pairs of points in $Z$.
   Then we have
   $$\alpha((2m-1)Z)=3m-1,\; \alpha(2mZ)=3m,\; \alpha((2m+1)Z)=X,\; \alpha(2(m+1)Z)=3(m+1).$$
   This show that $X$ is either $3m+1$ or $3m+2$. The first case is excluded
   by Corollary \ref{cor:alphas of components} as the points are not collinear.
   The remaining case is our claim.
\endproof
\section{Differences equal $2$}
   In this section we study subschemes $Z$ with differences equal to $2$.
   We begin by the following general statement, which restricts possible
   configurations. It parallels the result of Theorem \ref{thm:only lines}
   in the present setting.
\begin{theorem}\label{thm:rational components}
   Let $Z=\left\{P_1,\dots,P_r\right\}\subset\P^2$ be a finite set of points such that
   $$\alpha_{k,k-1}(Z)=2$$
   for some $k\geq 2$. Moreover, assume that $C$ is a divisor of degree
   $d=\alpha(kZ)$ vanishing along $kZ$. Then every irreducible component of $C$ is a rational curve;
\end{theorem}
\proof
  Let $\widetilde{C}$ be a component of $C= D + \widetilde{C}.$
  Let $\widetilde{d} = \deg \widetilde{C}$, $\widetilde{m_j} = \mult_{P_j} \widetilde{C}$.
  We want to show that $\widetilde{C}$ is rational.
  If $\widetilde{d}=1$ or $2$ then we are done. Similarly, if $\widetilde{d}=3$
  and $\widetilde{C}$ is singular, then we are also done. If $\widetilde{d}=3$
  and $C$ is smooth, then $D$ is a curve of degree $d-3$ vanishing
  along $Z$ to order $(k-1)$, which contradicts our assumption.
  Thus, we can assume that $\widetilde{d}\geq 4$.

   Lemma \ref{lem:alphas on subdivisors} implies that $h^0(\ldthreem)=0$.
   Counting conditions, we have
  $$(\widetilde{d}-3)\widetilde{d} < \sum_{j=1}^{r} (\widetilde{m_j}-1)\widetilde{m_j}.$$
  On the other hand, $\widetilde{C}$ is irreducible, hence by the genus formula
  $$(\widetilde{d}-1)(\widetilde{d}-2) \geq \sum_{j=1}^{r} (\widetilde{m_j}-1)\widetilde{m_j}.$$
  Observe that the above implies
  $$(\widetilde{d}-1)(\widetilde{d}-2)-1 \leq \sum_{j=1}^{r} (\widetilde{m_j}-1)\widetilde{m_j} \leq (\widetilde{d}-1)(\widetilde{d}-2),$$
  and since the terms in the middle and on the right hand side are always even, while the term on the left is odd, we have the equality in the
  genus formula. This implies that $\widetilde{C}$ has geometric genus $0$.
  Moreover, all singularities of $\widetilde{C}$ are ordinary multiple points.
\endproof
   Some examples of subschemes $Z$ with $\alpha_{2,1}=2$ have been discussed
   in \cite{BocCha11}. We repeat here some of them and discuss a few new ones.
   The first one must be expected after the statement of Corollary \ref{cor:collinear}.
\begin{example}\label{exp:conics}
   Let $Z=\left\{P_1,\dots,P_r\right\}$ be a set of $r\geq 5$ points on a smooth conic $C$.
   Then
   $$\alpha(kZ)=2k$$
   for all positive integers $k$.
\end{example}
   Heading for examples of irreducible curves satisfying Theorem \ref{thm:rational components}
   we have the following.
\begin{example}\label{ex:irreducible high sing}
   Let $d \geq 2$, $k \geq 2$, $r \geq 1$ be integers satisfying
   $$(d-1)(d-2)=rk(k-1), \qquad 2(d-1)<rk.$$
   Assume that there exists an irreducible curve $C \in H^0(\ldk)$. If $Z$ is the set of
   singular points of $C$ then
   $$\alpha_{k,k-1}(Z)=2.$$
\end{example}
\proof
   We show that $\alpha(kZ)=d$ and $\alpha((k-1)Z)=d-2$. To this end it suffices to show that
  \begin{eqnarray*}
  H^0(\P^2,\calo_{\P^2}(d)\otimes\cali(kZ)) \neq \emptyset , \\
  H^0(\P^2,\calo_{\P^2}(d-1)\otimes\cali(kZ)) = \emptyset , \\
  H^0(\P^2,\calo_{\P^2}(d-2)\otimes\cali((k-1)Z)) \neq \emptyset , \\
  H^0(\P^2,\calo_{\P^2}(d-3)\otimes\cali((k-1)Z))) = \emptyset .
  \end{eqnarray*}
  The first claim is satisfied by assumption.
  To prove the second, assume that there exists a curve $D$ of degree $d-1$ vanishing at $kZ$. Then, since $C$ is irreducible,
  by Bezout we have
  $$d(d-1) \geq rk^2 = rk(k-1)+rk > (d-1)(d-2)+2(d-2)=(d+1)(d-2),$$
  which is false.
  To prove the third claim, compute conditions
  $$(d-1)(d+2)-rk(k-1)=(d-1)(d+2)-(d-1)(d-2) > 0.$$
  To prove the last claim, again use Bezout
  $$d(d-3) \geq rk(k-1) = (d-1)(d-2),$$
  which gives absurd.
\endproof
   Examples \ref{exp:k2} and \ref{exp:k3} below show that curves satisfying
   numerical conditions of Example \ref{ex:irreducible high sing} exist for $k=2$
   and $k=3$.
   We do not  know if such curves exist
   for $k\geq 4$.
   Possibly examples to consider would be
   \begin{itemize}
      \item a curve of degree $13$ with $11$ points of multiplicity $4$;
      \item a curve of degree $17$ with $12$ points of multiplicity $5$;
      \item a curve of degree $22$ with $14$ points of multiplicity $6$ and so on.
   \end{itemize}
\begin{example}\label{exp:k2}
   Let $C$ be an irreducible curve of degree $d\geq 5$ with the maximal possible number of nodes,
   i.e. $\frac{(d-1)(d-2)}{2}$ of them. Such a curve exists by a classical result of Severi,
   see \cite{GLS98} for a modern approach and much more.
   Let $Z$ be the set of all nodes of $C$. Then
   $$\alpha(2Z)=d\;\; \mbox{ and }\;\; \alpha(Z)=d-2.$$
\end{example}
   The construction of an example with $k=3$ is a way more involved.
   We thank Joaquim Ro\'e for making us aware of papers \cite{Gia84} and \cite{GraMez88}
   and explaining to us their content.
\begin{example}\label{exp:k3}
   Let $C$ be a smooth plane cubic. There are $9$ inflection points on $C$
   and there are $12$ lines arranged so that each inflection point is contained
   in exactly $4$ lines and each line contains exactly $3$ points. This is
   the Hesse configuration, see \cite{ArtDol09} for more details and
   interesting background.

   The dual configuration consists of $9$ lines $L_1,\dots,L_9$ and $12$ points
   arranged so that each line passes through exactly $4$ points
   and each point is contained in exactly $3$ lines.
   The divisor $D'=L_1+\dots+L_9$ has degree $9$ and exactly $12$ triple points.
   Let $L$ be an additional line intersecting $D'$ transversally,
   in particular not passing through any of the triple points.
   The divisor $D=D'+L$ has degree $10$ and contains $12$ triple
   points and $9$ double points $L_i\cap L$. All of these singularities
   are ordinary, as they arise as intersection points of reduced lines.
   Let $W=W(10;3^{\times 12},2^{\times 9})$ be the variety parameterizing
   all plane curves of degree $10$ with $12$ triple points
   and $9$ double points. Let $V\subset W$ be a component containing $D$.
   Then \cite[Proposition 1.2]{Gia84} implies that $V$ has a good dimension
   (for precise definition we refer to \cite{Gia84} and \cite{GraMez88}).

   The divisor $D$ considered as an element of $W(10;3^{\times 12})$
   is virtually connected i.e. remains connected after removing arbitrary
   number of triple points. This is guaranteed by the line $L$ intersecting
   all components of $D'$. Then \cite[Theorem 4.3]{GraMez88} implies
   that a general member of $W(10; 3^{\times 12})$ is irreducible. Let
   $C$ be such a member and let $Z$ be the set of triple points on $C$.
   Then
   $$\alpha(2Z)=8\;\;\mbox{ and }\; \alpha(3Z)=10.$$
\end{example}
   We can derive new examples combining curves discussed in Example \ref{exp:k2}.
   For the purpose
   of this example we allow a smooth conic as a nodal rational curve (with
   no nodes). A line however is excluded.
\begin{example}
   Let $C_1, C_2$ be two nodal rational curves of degree $d_1\geq 2$ and $d_2\geq 2$ respectively,
   intersecting transversally. Let
   $$Z=\sing(C_1)\cup\sing(C_2)\cup(C_1\cap C_2).$$
   Then
   $$\alpha(2Z)=d_1+d_2\;\;\mbox{ and }\;\; \alpha(Z)=d_1+d_2-2.$$
\end{example}
\proof
   Assume that $\alpha(2Z)=\gamma\leq d_1+d_2-1$ and let $\Gamma$ be a divisor in
   $|\calo_{\P^2}(\gamma)\otimes \cali(2Z)|$. We show to the contradiction
   that $C_1$ and $C_2$ are components of $\Gamma$. Indeed, if not, then we have
   $$(d_1+d_2-1)d_1\geq \Gamma\cdot C_1\geq 2(d_1-1)(d_1-2)+2d_1d_2.$$
   This is equivalent to
   $$d_1(d_1+d_2-5)+4\leq 0,$$
   which is never satisfied for $d_1,d_2\geq 2$. The same argument applies to $C_2$.\\
   The second assertion is proved in the similar way.
\endproof

   Finally, we show that the irreducible curves appearing in Theorem \ref{thm:rational components}
   can all be lines. This implies in particular, that the converse to Theorem \ref{thm:only lines} b)
   cannot be true.
\begin{example}
   Let $L_1,\dots,L_d$ be a general arrangement of $d\geq 3$ lines in $P_2$ and let
   $Z$ be the set of points consisting of all intersection points $P_{ij}=L_i\cap L_j$
   but $P_{12}$. Then
   $$\alpha(2Z)=d\;\; \mbox{ and }\; \alpha(Z)=d-2.$$
\end{example}
\proof
   The proofs of both claims are the same. If there were a divisor of a lower
   degree in either case, then we show, intersecting it with configuration
   lines, that it must contain all lines in the first case and the lines $L_3,\dots,L_d$
   in the second case.
\endproof
   The next example shows that there is no straightforward generalization of
   Corollary \ref{cor:collinear}, there might be two consecutive differences of initial degree
   equal to $2$ without the points being forced to lie on a conic. Lemma \ref{lem:after}
   exhibits another example of this kind with $3$ consecutive differences equal $2$,
   but the argument there is less explicit.
\begin{example}\label{ex:two jumps by two}
   Let $Z=\left\{P_1,\dots,P_6\right\}$ be a set of $6$ general points
   in the plane. Let $C_i$ be the conic passing through all points in $Z$ but $P_i$,
   for $i=1,\dots,6$. Then
   $$\alpha(5Z)=12,\; \alpha(4Z)=10\; \mbox{ and }\; \alpha(3Z)=8.$$
\end{example}
\proof
   Again, the argument is the Bezout's Theorem. We discuss only the last case.
   Assume to the contrary that there is a divisor $\Gamma_1$ of degree at most $7$ vanishing
   along $3Z$. If $C_1$ is not a component of $\Gamma_1$, then intersecting with $C_1$ we have
   $$14\geq \Gamma_1\cdot C_1\geq 5\cdot 3=15,$$
   a contradiction. Thus $C_1$ is contained in $\Gamma_1$ and there is a divisor
   $\Gamma_2=\Gamma_1-C_1$ of degree at most $5$ with vanishing vector $(3,2,2,2,2,2)$ at $P_1,\dots,P_6$.
   If $C_2$ is not a component of $\Gamma_2$, then we have
   $$10\geq \Gamma_2\cdot C_2\geq 3+4\cdot 2=11,$$
   a contradiction. Thus $C_2$ is a component of $\Gamma_2$ and we get a new divisor
   $\Gamma_3=\Gamma_2-C_2$ of degree at most $3$ and vanishing vector $(2,2,1,1,1,1)$.
   This divisor must split the line $L_{12}$ through the points $P_1$ and $P_2$.
   Then the residual curve $\Gamma_4=\Gamma_3-L_{12}$ is a conic containing all
   the points $P_1,\dots,P_6$, which contradicts our assumption that the points
   are general.

   Note, that in this example $\alpha(2Z)\leq 5$. Indeed, it is enough to take
   the conics $C_1$ and $C_2$ and the line $L_{12}$ to obtain a divisor
   with double points along $Z$.
\endproof
   It is nevertheless natural to wonder if there is a result along the lines
   of Corollary \ref{cor:collinear} in the case of degree jumping by $2$.
\begin{conjecture}
   Assume that $k\geq 5$ and $Z$ is a collection of points such that
   $$\alpha_{k,k-1}(Z)=\alpha_{k-1,k-2}(Z)=\alpha_{k-2,k-3}(Z)=\alpha_{k-3,k-4}(Z)=2.$$
   Then $Z$ is contained in a single conic.
\end{conjecture}
   We were not able to prove this conjecture, yet there is a strong supporting
   evidence, steaming partly from the next lemmata and Theorem \ref{thm:last thm}.

   Now, we want to investigate more closely rational curves appearing
   as components in Theorem \ref{thm:rational components}. Taking into account
   Corollary \ref{cor:alphas of components}, we arrive at the following statement.
\begin{lemma}\label{lem:excess curve}
   Let $C$ be an irreducible curve of degree $d$ with multiplicities
   $m_1\geq 2$, $m_2,\dots,m_r$ at points $P_1,\dots,P_r$. Let
   $\subn=(m_1-2,m_2-1,\dots,m_r-1)$. Then
   $$\alpha(\subn Z)\leq d-3,$$
   where as usually $Z=\left\{P_1,\dots,P_r\right\}$.
\end{lemma}
\proof
   It follows from the genus formula that
   $$(d-1)(d-2)\geq \sum_{i=1}^rm_i(m_i-1).$$
   The existence of a divisor of degree $d-3$ with vanishing vector $\subn$
   follows from the inequalities
   $$\begin{array}{rcccl}
      (m_1-2)(m_1-1)+\sum\limits_{i=2}^rm_i(m_i-1) & = & \sum\limits_{i=1}^r m_i(m_i-1)-2(m_1-1) & \leq & \\
      \\
      & \leq & d(d-3)+2(2-m_1) & \leq & d(d-3)
   \end{array}$$
   and we are done.
\endproof
   As a consequence we derive the following statement.
\begin{lemma}\label{lem:excess divisor}
   Let $Z=\left\{P_1,\dots,P_r\right\}$ be a reduced $0$--dimensional
   subscheme of $\P^2$. Let $C=\sum a_jC_j$ be a divisor of degree $d=\alpha(kZ)$
   vanishing at $kZ$ for some $k\geq 2$. Assume that $\alpha((k-1)Z)=d-2$ and
   that for each point $P_i$ there exists at least one component $C_{j(i)}$
   singular at $P_i$, i.e. $\ord_{P_i}C_{j(i)}\geq 2$. Then
   $$\ord_{P_i}C=k$$
   for all $i=1,\dots,r$.
\end{lemma}
\proof
   Suppose to the contrary that, after possible renumbering of the points, $\ord_{P_1}C\geq k+1$
   and, after possible renumbering of the components, $\ord_{P_1}C_1\geq 2$. Let
   $m_i=\ord_{P_i}C_1$. For the divisor $D=C-C_1$ we have
   $$\ord_{P_i}D\geq k-m_i\;\;\mbox{ for }\;\;i\geq 2\;\;\mbox{ and }\;\; \ord_{P_1}D\geq k+1-m_1.$$
   Lemma \ref{lem:excess curve} applies to the curve $C_1$, so that
   there exists a divisor $\Gamma$ of degree  $\deg(\Gamma)=\deg(C_1)-3$ with
   $$\ord_{P_i}\Gamma\geq m_i-1\;\;\mbox{ for }\;\;i\geq 2\;\;\mbox{ and }\;\; \ord_{P_1}\Gamma\geq m_1-2.$$
   Then $\Gamma+D$ has degree $d-3$ and vanishes along $(k-1)Z$, which contradicts our assumptions.
\endproof
   A straightforward useful consequence of the above two lemmata is the following.
\begin{corollary}\label{cor:excess}
   Let $Z=\left\{P_1,\dots,P_r\right\}$ be a reduced $0$--dimensional
   subscheme of $\P^2$. Let $C=\sum a_jC_j$ be a divisor of degree $d=\alpha(kZ)$
   vanishing at $kZ$ for some $k\geq 2$. Assume that for a fixed $i\in\left\{1,\dots,r\right\}$
   there is a component $C_i$ of $C$ singular at $P_i$ and that $\ord_{P_i}C\geq k+1$. Then
   $\alpha((k-1)Z)\leq d-3$.
\end{corollary}
   Now we are heading for a result paralleling Theorem \ref{thm:first thm}
   in the present setting. First we need some preparations. The following
   Lemma and its proof are motivated by a result of Lazarsfeld \cite[Proposition 2.5]{Laz10}.
\begin{lemma}\label{lem:Zariski like}
   Let $D\subset \P^2$ be a reduced divisor of degree $d$ with multiplicities
   $m_1\geq\dots\geq m_r\geq 2$ in some points $P_1,\dots,P_r$. Let
   $\subm=(m_1,\dots,m_r)$ be the multiplicities vector and $Z=\left\{P_1,\dots,P_r\right\}$.
   Then $((\subm-\one)Z)$ imposes independent conditions on curves of degree $k\geq d-2$,
   i.e.
   $$H^1(\P^2,\calo_{\P^2}(k)\otimes \cali((\subm-\one)Z))=0$$
   for all $k\geq (d-2)$.
\end{lemma}
\proof
   Let $\Gamma$ be a reduced divisor
   of degree $\gamma$ vanishing along $Z$ with no common components with $D$.
   For the $\Q$--divisor
   $$D'=(1-\varepsilon)D+m_1\varepsilon \Gamma,$$
   we have
   $$\ord_{P_i}D'=(1-\varepsilon)\ord_{P_i}D + m_1\varepsilon\ord_{P_i}\Gamma\geq m_i$$
   and
   $$\deg(D')=(1-\varepsilon)d+m_1\varepsilon \gamma<d+1$$
   for $\varepsilon$ sufficiently small.
   For the line bundle $L=\calo_{\P^2}(k+3)$, the difference $(L-D')$ is big and nef
   for all $k\geq d-2$
   and we have
   $$H^1(\P^2, \calo_{\P^2}(k)\otimes\cali(D'))=0$$
   by Nadel vanishing theorem \cite[Theorem 9.4.8]{PAG}.
   On the other hand, the multiplier ideal $\cali(D')$ vanishes on $(\subm-\one)Z$
   because of \cite[Proposition 9.3.2]{PAG} and it is cosupported
   on a finite set since no component of $D'$ has a coefficient $\geq 1$.
\endproof
\begin{theorem}\label{thm:last thm}
   Let $Z\subset\P^2$ be a finite set of $r$ points and let $k\geq 5$ be an integer.
   If $p$ is an integer such that
   $$\alpha(Z)=p,\;\; \alpha(2Z)=p+2,\;\;\dots,\alpha(kZ)=p+2(k-1),$$
   then $p=2$, i.e. the points in $Z$ are contained in a single conic.
\end{theorem}
\proof
   If $r\leq 5$, then there is nothing to prove.
   We postpone the case $r=6$ until Lemma \ref{lem:after}. So we assume from now on
   \begin{equation}\label{eq:r geq 7}
   r\geq 7.
   \end{equation}
   It follows from the proof below that in this case one can relax the assumption on $k$ to $k\geq 4$.

   Let $C_i$ denote a divisor of minimal degree vanishing along $iZ$ for $i=1,\dots,k$.
   From \cite[Proposition 4.4]{BocCha11} we know that for $C_2$ there are the
   following possibilities:
   \begin{itemize}
      \item[a)] $C_2$ is a double conic and $p+2=4$, or
      \item[b)] $C_2$ contains a double line $2L$ and the set
      $Z'=Z\setminus L$ is a non-empty configuration of all
      intersection points of $p$ general lines $L_1,\dots,L_p$, or
      \item[c)] $C_2$ is reduced.
   \end{itemize}
   In case a) we are immediately done.\\
   In case b) we observe first that there are $s\geq p$ points from $Z$ on the line $L$.
   Indeed, as there is no curve of degree $(p-1)$ containing $Z$, it must be
   $$0\geq {{p+1}\choose{2}}-{{p}\choose{2}}-s=p-s.$$
   Next, let $D=C_4$ be a divisor of degree $p+6$ vanishing along $4Z$. We have
   $$p+6=D\cdot L\geq 4\cdot s\geq 4\cdot p,$$
   which implies either that $p=2$ in which case we are done, or that $L$ is a component of $D$.
   Similarly, step by step, we show that $L$ is a component of $(D-L)$ and that
   $(D-2L)$ contains all lines from the arrangement. Thus
   $$D-2L-L_1-\dots-L_p$$
   is a divisor of degree $4$ vanishing along $2Z$. It follows again, that
   $p=2$ and we are done.

   \begin{center}
      \emph{For the rest of the proof we assume that $C_2$ is a reduced divisor.}
   \end{center}
   Then Lemma \ref{lem:Zariski like} implies that the points $P_1,\dots,P_r$
   impose independent conditions on curves of degree $p$. Hence
   at most $(p+1)$ of them lie on a line and we have the following
   bounds on their number
   \begin{equation}\label{eq:bound on r}
      {{p+1}\choose{2}}\leq r\leq {{p+2}\choose{2}}-1.
   \end{equation}
   \begin{center}
   \textbf{Case $C_3$ reduced.}
   \end{center}
   If $C_3$ is reduced, we derive from
   Lemma \ref{lem:Zariski like} and the existence of $C_2$ that
   \begin{equation}\label{eq:2bound on r}
      {{p+4}\choose{2}}-3r\geq 1,
   \end{equation}
   which, taking into account \eqnref{eq:bound on r}, implies $p\leq 3$.
   If $p\leq 2$, we are done. If $p=3$, then \eqnref{eq:2bound on r}
   implies $r\leq 6$, which contradicts \eqnref{eq:r geq 7}.
   \begin{center}
   \textbf{Case $C_3$ non-reduced.}
   \end{center}
   Let
   $$C_3=a\Gamma+R$$
   be a decomposition of $C_3$ with $a=2$ or $a=3$ (note
   that $a>3$ is immediately excluded as then one could replace
   $C_3$ by a lower degree divisor $(a-1)\Gamma+R$). We
   assume that $\Gamma$ and $R$ have no common components. We write
   $Z'=Z\setminus\Gamma$ and $Z''=Z\setminus Z'$, so that
   $Z'$ is contained entirely in $R$ and $Z''$ is contained in $\Gamma$.
   Corollary \ref{cor:excess} implies that $\Gamma$ is smooth in
   points of $Z$. Hence $(a-1)\Gamma+R$ vanishes along $2Z$. This shows
   that $\Gamma$ is either a line or a (possibly singular) conic.
   \begin{center}
   \textbf{Subcase $C_3$ contains a triple component, i.e. $a=3$.}
   \end{center}
   If $R=0$, then we are done. Otherwise $R$ is a divisor of degree
   $p+4-3\deg(\Gamma)$ vanishing to order $3$ along $Z'$.
   Taking derivatives of the equation of $R$, we see that there
   is another divisor $R''$ of degree $p+2-3\deg(\Gamma)$ vanishing
   along $Z'$. The union $R''\cup \Gamma$ has degree $p+2-2\deg(\Gamma)$
   and vanishes along $Z$. This shows that $\deg(\Gamma)=1$ in this case.
   Then the divisor $R$ has degree $p+1$ and vanishes
   to order $3$ on $Z'$. Since there is no
   divisor of degree $\leq p-2$ vanishing on $Z'$ (otherwise its union with $\Gamma$
   would have degree $p-1$ and would vanish on $Z$), we have
   $$\alpha(Z')=p-1,\;\; \alpha(2Z')=p,\;\; \alpha(3Z')=p+1.$$
   Theorem \ref{thm:first thm} implies that $Z'$ is contained in a line, hence $Z$ is contained in a conic.
   \begin{center}
   \textbf{Subcase $C_3$ contains a double component, i.e. $a=2$.}
   \end{center}
   We study now the cases $\deg(\Gamma)=1$ and $\deg(\Gamma)=2$ separately.

   \textbf{Subsubcase $\Gamma$ is a conic.}

   In this case we have $\deg(R)=p$ and $R$ has multiplicity $3$ along $Z'$.
   Taking twice a derivative of the equation
   of $R$, we see that there exists a curve of degree $p-2$ vanishing along $Z'$
   On the other hand, there cannot exist a curve of lower degree, as its union with
   $\Gamma$ would be of degree less than $p$ and would vanish along $Z$.
   Thus Corollary \ref{cor:collinear} implies
   that $Z'$ consists of collinear points. Let $L$ be a line containing $Z'$.
   Then $\Gamma+L$ vanishes along $Z$. Hence $p=3$.

   The divisor $R$ has degree $3$, vanishes with multiplicity $3$ along $Z'$
   and it vanishes along $Z''$. If there are at least $2$ points
   in the set $Z'$, then $R=3L$ and it cannot vanish along $Z''$
   (if $Z''$ is empty, then $Z=Z'$ is contained in a line and we are done). Hence $Z'$
   is a single point. There must be at least $6$ points in $Z''$, as there
   are at least $7$ points altogether. On the other hand, the divisor $R$
   consists of $3$ lines $L_1, L_2, L_3$ whose union vanishes in all points in $Z''$. This shows
   that there are at most $6$ points in $Z''\subset \Gamma$.

   So there are exactly $6$ points from $Z$ on the conic $\Gamma$, the intersection points
   of $\Gamma$ with the lines $L_1,L_2,L_3$.
   We show that then $\alpha(4Z)\geq 10$, contradicting the assumptions
   of the theorem.
   Suppose to the contrary that there exists a divisor $D$ of degree $9$
   vanishing along $4Z$. Then, by the standard Bezout argument, we show
   step by step that $D$ must split off the curves $2\Gamma$, $L_1$, $L_2$ and $L_3$.
   The residual divisor has than degree $2$ and it vanishes in all points of $Z$,
   a contradiction.

   \textbf{Subsubcase $\Gamma$ is a line.}

   In this situation we have $\deg(R)=p+2$ and $R$ has multiplicities $3$ in points
   of the set $Z'$. Note that there is no divisor of degree
   less or equal $p-2$ vanishing along $Z'$, since its union with $\Gamma$
   would be of degree $p-1$ and would vanish along $Z$. This gives the following
   possibilities for $\alpha$'s of $Z'$:
   \begin{center}
   \begin{tabular}{|l||c|c|c|}
      \hline
      Case & $\alpha(Z')$ & $\alpha(2Z')$ & $\alpha(3Z')$ \\
      \hline
      \hline
      A) & $p$ & $p+1$ & $p+2$\\
      B) & $p-1$ & $p$ & $p+1$\\
      C) & $p-1$ & $p$ & $p+2$\\
      D) & $p-1$ & $p+1$ & $p+2$\\
      \hline
   \end{tabular}
   \end{center}
   In cases A) and B), Theorem \ref{thm:first thm} implies that $Z'$
   is contained in a line, hence $Z$ is contained in a union of two
   lines and we are done.

   In case C), the main result in \cite{BocCha11} implies that
   either the points in $Z'$ are all collinear and we are done,
   or that
   $Z'$ consists of all intersection points of a general arrangement
   of $p$ lines $L_1,\dots,L_p$. In the later case, there are exactly
   ${{p}\choose{2}}$ points in $Z'$. Hence by \eqnref{eq:bound on r}
   there are at least
   $${{p+1}\choose{2}}-{{p}\choose{2}}=p$$
   points in $Z''$. The divisor $D=\Gamma+R$ is reduced and
   has multiplicities $2$ on $Z''$ and multiplicities $3$ on $Z'$.
   Lemma \ref{lem:Zariski like} implies that $Z''+2Z'$ imposes independent
   conditions on curves of degree $(p+1)$. Since
   $$\Gamma\cup L_1\cup\dots\cup L_p$$
   is a divisor of degree $p+1$ vanishing on $Z''+2Z'$, it must be
   $${{p+3}\choose{2}}-p-3{{p}\choose{2}}\geq 1.$$
   This is possible only if $p=3$, so that there are exactly $3$ points in
   $Z'$ and at least 4 points on $\Gamma$. Let $L_1$, $L_2$, $L_3$ be the lines
   determined by the pairs of points in $Z'$. Consider $C_4$ of degree 9 passing through $4Z$.
   It follows (by repeated use of Bezout) that $C_4=3\Gamma+2L_1+2L_2+2L_3$,
   but it is impossible, since then at most three points on $\Gamma$ can attain the multiplicity 4.

    In case D) $R$ splits into $(p+2)$ lines by Theorem \ref{thm:only lines}.
    None of these lines is a double line, as this case was already covered by b2).
    This arrangement of lines has triple points in $Z'$ and vanishes in $Z''$.
    It follows that there are at most $(p+2)$ points in $Z''$ (intersection points
    of $\Gamma$ with the arrangement lines) and at most $\frac13{{p+2}\choose{2}}$
    points in $Z'$ (we use here a very rough estimate on the number of possible
    triple points in a configuration of lines). Taking into account \eqnref{eq:bound on r},
    we must have
    $$p+2+\frac13{{p+2}\choose{2}}\geq{{p+1}\choose{2}}.$$
    This is possible only for $p=3$ or $p=4$. If $p=3$, then there are
    at most $2$ points in $Z'$, hence $\alpha(Z')=1=p-2$, a contradiction.
    Similarly, if $p=4$, then there are at most $3$ points in the set $Z'$
    and $\alpha(Z')=2=p-2$ again.
\endproof
   Now we turn to the case of $6$ points. This Lemma shows in particular that
   the assumptions in Theorem \ref{thm:last thm} are sharp.
\begin{lemma}\label{lem:after}
   If $Z$ consists of 6 points and
   $$\alpha(Z)=p,\; \alpha(2Z)=p+2,\; \alpha(3Z)=p+4,\; \alpha(4Z)=p+6,$$
   then either $p=2$ and $Z$ lies on a conic or $p=3$ and $Z$ is a configuration of
   intersection points $A$, $B$, $C$ of three general lines,
   and additional points $D$, $E$, $F$, each of those lying on exactly one of the lines,
   one for each line. This is indicated on the picture below
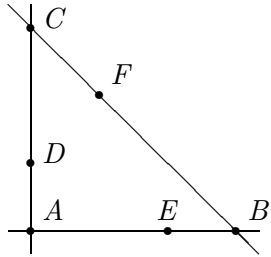
\unitlength.15mm
\begin{center}
\begin{picture}(220,220)(0,0)
\put(20,0){\line(0,1){220}}   
\put(0,20){\line(1,0){220}}   
\put(0,220){\line(1,-1){220}}   
 \put(20, 200){\circle*{7}} \put(30,200){C}
 \put(20, 80){\circle*{7}} \put(30,80){D}
 \put(20, 20){\circle*{7}} \put(30,30){A}
 \put(140, 20){\circle*{7}} \put(130,30){E}
 \put(200, 20){\circle*{7}} \put(210,30){B}
 \put(80, 140){\circle*{7}} \put(90,150){F}
\end{picture}
\end{center}
   If moreover
   $$\alpha(5Z)=p+8,$$
   then $Z$ is contained in a conic, i.e. $p=2$.
\end{lemma}
\proof
   Since 6 points always lie on a cubic, we have $p \leq 3$ and the only interesting case is $p=3$.

   We say that two sets $Z=\left\{P_1,\dots,P_6\right\}$, $Z'=\left\{P_1',\dots,P_6'\right\}$,
   of points in $\P^2$ are equivalent if they have the same Hilbert functions for any set of multiplicities, i.e
   $$ \dim (I(m_1P1+\ldots+m_6P_6))_t = \dim (I(m_1P_1'+\ldots+m_6P_6'))_t$$
   for all non-negative integers $t$, $m_1,\dots,m_6$.
   An equivalence class is called a type.

   Obviously, equivalent sets of points have the same initial degree for any multiplicity sequence.

   From \cite{GuaHar07} we know that for $6$ points there are exactly $11$ types.
   Moreover, there are algorithms to find Hilbert functions for a given type and multiplicities
   \cite{HarAlg}. Using them we verify that in order to have $\alpha(4Z)=9$,
   the only possibility is to have $Z$ of type $9$ in \cite{HarAlg},
   which is equivalent to the configuration described in the Lemma.
   The picture below gives examples of divisors vanishing along $2Z$, $3Z$ and $4Z$.
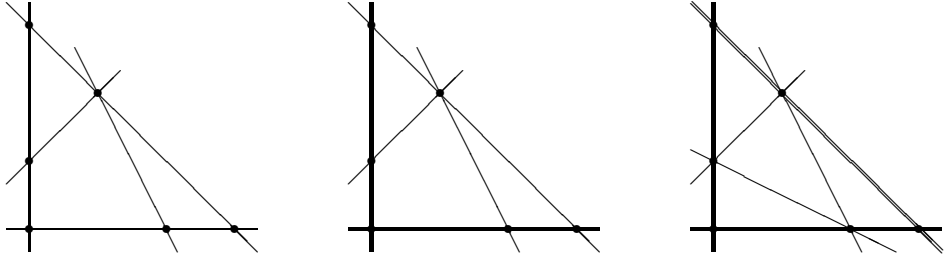
\begin{center}
\unitlength.15mm
\begin{picture}(780,220)(0,0)
\put(20,0){\line(0,1){220}}   
\put(0,20){\line(1,0){220}}   
\put(0,220){\line(1,-1){220}}   
\put(60,180){\line(1,-2){90}}   
\put(0,60){\line(1,1){100}}   
 \put(20, 200){\circle*{7}}
 \put(20, 80){\circle*{7}}
 \put(20, 20){\circle*{7}}
 \put(140, 20){\circle*{7}}
 \put(200, 20){\circle*{7}}
 \put(80, 140){\circle*{7}}
\put(321,0){\line(0,1){220}}   
\put(319,0){\line(0,1){220}}   
\put(300,21){\line(1,0){220}}   
\put(300,19){\line(1,0){220}}   
\put(300,220){\line(1,-1){220}}   
\put(360,180){\line(1,-2){90}}   
\put(300,60){\line(1,1){100}}   
 \put(320, 200){\circle*{7}}
 \put(320, 80){\circle*{7}}
 \put(320, 20){\circle*{7}}
 \put(440, 20){\circle*{7}}
 \put(500, 20){\circle*{7}}
 \put(380, 140){\circle*{7}}
\put(621,0){\line(0,1){220}}   
\put(619,0){\line(0,1){220}}   
\put(600,21){\line(1,0){220}}   
\put(600,19){\line(1,0){220}}   
\put(600,219){\line(1,-1){220}}   
\put(601,221){\line(1,-1){220}}   
\put(600,90){\line(2,-1){180}}   
\put(660,180){\line(1,-2){90}}   
\put(600,60){\line(1,1){100}}   
 \put(620, 200){\circle*{7}}
 \put(620, 80){\circle*{7}}
 \put(620, 20){\circle*{7}}
 \put(740, 20){\circle*{7}}
 \put(800, 20){\circle*{7}}
 \put(680, 140){\circle*{7}}
\end{picture}
\end{center}
   Further computations using the above mentioned algorithms show that
   for this type $\alpha(5Z)=12$. This concludes the proof of the Lemma and thus
   also the proof of Theorem \ref{thm:last thm}.
\endproof
\begin{remark}
   Note that one could prove Lemma \ref{lem:after} along the lines of Example \ref{ex:two jumps by two}
   without calling to a computer aided argument. However, this would take several pages
   of dull computations which we preferred to avoid.
\end{remark}
   Now we show that one cannot expect a result along lines of Theorem stated
   in the Introduction for $d\geq 4$.
\begin{example}\label{ex:Nagata}
   Let $Z$ be a set of $16$ very general points in the projective plane.
   Then \cite[Proposition]{Nag59} implies that for all $k\geq 1$
   there is no divisor of degree $4k$ vanishing to order $k$ along $kZ$.
   On the other hand, an elementary conditions counting shows that there
   is always a divisor of degree $4k+1$ vanishing along $kZ$. These two
   facts together imply that
   $$\alpha_{k+1,k}=4$$
   for all $k\geq 1$ but $Z$ is not contained in a curve of degree $4$
   (since $\alpha(Z)=5$).
\end{example}
   With a little more care the same idea can be used to produce similar examples
   for arbitrary $d\geq 4$.

\paragraph*{Acknowledgement.} We would like to thank Thomas Bauer, Joaquim Ro\'e and
   Stefan M\"uller-Stach for helpful remarks and comments. Parts of
    this paper
   were written while the second author was a visiting professor at the University
   Mainz as a member of the program Schwerpunkt Polen. The nice working
   conditions and financial support are kindly acknowledged.


\bigskip \small

\bigskip
   Marcin Dumnicki,
   Jagiellonian University, Institute of Mathematics, {\L}ojasiewicza 6, PL-30348 Krak\'ow, Poland

\nopagebreak
   \textit{E-mail address:} \texttt{Marcin.Dumnicki@im.uj.edu.pl}

\bigskip
   Tomasz Szemberg,
   Instytut Matematyki UP,
   Podchor\c a\.zych 2,
   PL-30-084 Krak\'ow, Poland.

\nopagebreak
   \textit{E-mail address:} \texttt{szemberg@up.krakow.pl}

\bigskip
   Halszka Tutaj-Gasi\'nska,
   Jagiellonian University, Institute of Mathematics, {\L}ojasiewicza 6, PL-30348 Krak\'ow, Poland

\nopagebreak
   \textit{E-mail address:} \texttt{Halszka.Tutaj@im.uj.edu.pl}


\end{document}